\documentclass[12pt]{article}
\usepackage{mathtools} 
\usepackage{amsmath, amssymb} 
\usepackage{amsthm}
\usepackage{enumerate}  % lets you do \begin{enumerate}[(a)] and such
\usepackage{url} % formats url's nicely using \url{}

\usepackage{tikz}
\usepackage{graphicx}
\usepackage{caption,subcaption}

\newcommand\cx{{\mathbb C}}% complexes

\newcommand\cA{{\mathcal A}}

\DeclarePairedDelimiter\abs{\lvert}{\rvert}%
\DeclarePairedDelimiter\norm{\lVert}{\rVert}%

\makeatletter
\let\oldabs\abs
\def\abs{\@ifstar{\oldabs}{\oldabs*}}
\let\oldnorm\norm
\def\norm{\@ifstar{\oldnorm}{\oldnorm*}}
\makeatother

\newcommand\comp[1]{{\mkern2mu\overline{\mkern-2mu#1}}}

\newcommand\seq[4]{#1_{#2},#1_{#3},\ldots,#1_{#4}}

\newtheoremstyle{plainsl}%
	{\topsep}
	{\topsep}
	{\slshape} % only non-default setting
	{}
	{\normalfont\bfseries}
	{.}
	{ }
	{}

\swapnumbers

{\theoremstyle{plainsl}
\newtheorem{theorem}{Theorem}[section]
\newtheorem{lemma}[theorem]{Lemma}
\newtheorem{corollary}[theorem]{Corollary}}
{\theoremstyle{remark}
}

\renewcommand\proof{\noindent\textsl{Proof. }}
\newcommand\sqr[2]{{\vbox{\hrule height.#2pt
    \hbox{\vrule width.#2pt height#1pt \kern#1pt
        \vrule width.#2pt}\hrule height.#2pt}}}
\renewcommand\qed{%
	\ifmmode\eqno\sqr53
	\else\nolinebreak\ \hfill\sqr53\medbreak\fi}

\DeclareMathOperator{\LD}{LD}
\DeclareMathOperator{\dist}{dist}
\newcommand\ip[2]{\langle#1,#2\rangle}
\newcommand\one{{\bf1}}

\usepackage{blkarray}

\title{Factoring  Discrete-Time Quantum Walks on Distance Regular Graphs into  Continuous-Time Quantum Walks}

\author{Hanmeng Zhan \thanks{Department of Mathematics and Statistics, York University, Toronto, ON, Canada \texttt{h3zhan@yorku.ca}}}

\begin{document}
\maketitle

\begin{abstract}
We consider a  discrete-time quantum walk, called the Grover walk, on a distance regular graph $X$. Given that $X$ has diameter $d$ and invertible adjacency matrix, we show that the square of the transition matrix of the Grover walk on $X$ is a product of at most $d$ commuting transition matrices of  continuous-time quantum walks, each on some distance digraph of the line digraph of $X$. We also obtain a similar factorization for any graph $X$ in a Bose Mesner algebra.
\end{abstract}



\section{Introduction}
There are two types of quantum walks: continuous-time and discrete-time. In a  continuous-time quantum walk, the evolution is usually determined by the adjacency matrix or the Laplacian matrix of the graph. Hence, certain properties of the walk, such as perfect state transfer, uniform mixing and fractional revival, translate into properties of the graph (see, for example, \cite{Godsil2012,Chan2013,Chan2018}). On the other hand, the transition matrix of a  discrete-time quantum walk is a product of two non-commuting sparse unitary matrices, whose sizes usually equal the number of arcs of the graph. Moreover, the choice of these sparse matrices is not unique. Such freedom leads to various models of  discrete-time quantum walks (see, for example, \cite{Aharonov2000,Szegedy2004,Kendon2003,Portugal2015}); however, the relation between graph properties and walk properties is less clear.

In this paper, we study a  discrete-time quantum walk formalized by Kendon \cite{Kendon2003}. The transition matrix is a product of two matrices: the arc-reversal matrix, which maps an arc $(a, b)$ to the arc $(b, a)$, and a Grover coin matrix, which sends an arc $(a, b)$ to certain linear combination of the outgoing arcs of vertex $a$. Over the past few years, investigation of this model on specific families of graphs has received an increased attention. For example, Konno et al \cite{Konno2019} characterized the positive support of the $m$-th power of the transition matrix for regular graphs with girth greater than $2m-1$, and Yoshie \cite{Yoshie2019} characterized certain distance regular graphs on which this walk is periodic. Following their notation, we will refer to this walk as the \textsl{Grover walk}. 

The graphs we consider in this paper are distance regular, or, more generally, lying in the Bose Mesner algebra of some association scheme. Let $U$ be the transition matrix of a Grover walk. We show that if our graph $X$ is distance regular of diameter $d$ with invertible adjacency matrix, then
\[U^2 = \exp(t_1 S_1) \exp(t_2 S_2) \cdots \exp(t_d S_d),\]
where $S_i$ is the skew-adjacency matrix of the $i$-th distance digraph of the line digraph of $X$. Moreover, $S_i$ and $S_j$ commute. A similar factorization holds if $X$ lies in  a Bose-Mesner algebra of some association scheme, not necessarily constructed from a distance regular graph. This reveals an interesting connection between  discrete-time quantum walks on graphs and  continuous-time quantum walks on oriented graphs.

\section{Grover Walks}
To gather more combinatorial insights into the Grover walk, we start this section with an alternative definition given by Godsil and Guo \cite{Godsil2010Guo}.

\begin{figure}[h]
	\centering
\begin{minipage}{0.5\textwidth}
\includegraphics[width=6cm]{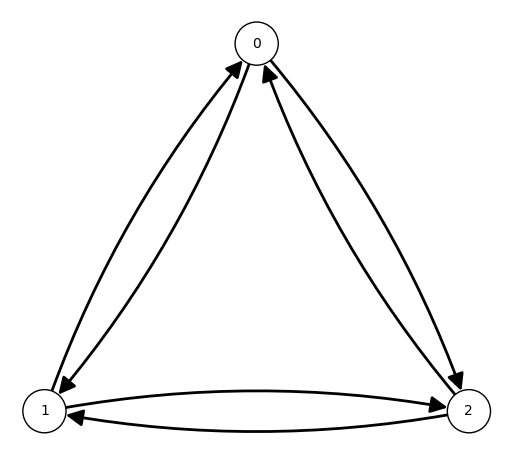}
\end{minipage}%
\begin{minipage}{0.5\textwidth}
\includegraphics[width=6cm]{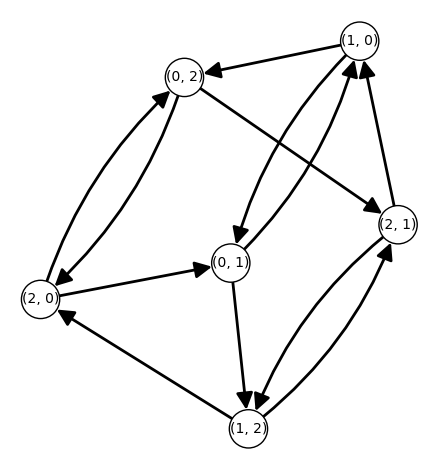}
\end{minipage}
\caption{$K_3$ and its line digraph}
\label{fig_LDK3}
\end{figure}

Throughout, let $X$ be a $k$-regular graph. We will replace each edge $\{a,b\}$ of $X$ by two arcs $(a,b)$ and $(b,a)$,   that is, the two directions of $\{a, b\}$. The \textsl{tail} and the \textsl{head} of an arc $(a,b)$ are $a$ and $b$, respectively. The \textsl{line digraph} of $X$, denoted $\LD(X)$, is the digraph whose vertices are the arcs of $X$, and $(a,b)$ is adjacent to $(c,d)$ in $\LD(X)$ if $b=c$. Note that the adjacency in $\LD(X)$ is not a symmetric relation. In Figure \ref{fig_LDK3}, we illustrate $K_3$ and its line digraph. 

Let $A(\LD(X))$ denote the  $(0,1)$-adjacency matrix of $\LD(X)$. Let $R$ denote the \textsl{arc-reversal matrix}, that is, the permutation matrix on the arcs that maps $(a,b)$ to $(b,a)$. We define the \textsl{transition matrix} of the Grover walk on $X$ to be
\begin{equation}\label{eq:defU}
U:= \frac{2}{k} A(\LD(X)) -R.
\end{equation}

While this definition assumes $X$ is regular, one can in fact define the Grover walk on any graph. Let $C$ be the matrix acting on $\cx^{\mathrm{arcs}(X)}$, such that for each arc $(u,v)$, its characteristic vector is sent by $C$ to a linear combination of the characteristic vectors of all arcs leaving $u$:
\[C e_{(u,v)} = \left(\frac{2}{\deg(u)}-1\right) e_{(u,v)} + \frac{2}{\deg(u)} \sum_{w\sim u, \, w\ne v} e_{(u,w)}.\]
This matrix $C$ is unitary and called the \textsl{coin matrix}. Given $R$ and $C$, the alternative definition of the transition matrix is simply $U=RC$ \cite{Kendon2003}. To see why this definition coincides with \eqref{eq:defU} when $X$ is regular, we introduce two more matrices, $D_t$ and $D_h$, whose rows are indexed by the vertices of $X$, and columns by the arcs of $X$:
\[(D_t)_{u, ab} = \begin{cases}
1,\quad u=a\\
0,\quad \text{otherwise}
\end{cases}\]
and
\[(D_h)_{u, ab} = \begin{cases}
1,\quad u=b\\
0,\quad \text{otherwise}
\end{cases}\]
We sometimes refer to $D_t$ as the \textsl{tail-arc} incidence matrix, and $D_h$ the \textsl{head-arc} incidence matrix. It is not hard to verity the following.

\begin{lemma}\label{lem_props}
Let $X$ be a $k$-regular graph. The incidence matrices $D_t$ and $D_h$ satisfy the following identities.
\begin{enumerate}[(i)]
\item $D_t R = D_h$
\item $D_tD_t^T = D_hD_h^T=kI$
\item $D_tD_h^T = D_hD_t^T = A(X)$
\item $D_h^TD_t = A(\LD(X))$
\qed
\end{enumerate}
\end{lemma}

Hence, 
\begin{align*}
U &= \frac{2}{k}A(\LD(X)) - R\\
&= \frac{2}{k} D_h^TD_t - R\\
&=R\left(\frac{2}{k}D_t^TD_t-I\right),
\end{align*}
where $\frac{2}{k}D_t^TD_t-I$ is precisely the coin matrix.

\section{Distances in digraphs}

Given a digraph $Y$,  the \textsl{distance} from $a$ to $b$, denoted $\dist_Y(a,b)$, is the length of the shortest dipath in $Y$ from $a$ to $b$. Note that in general $\dist_Y(a,b)\ne \dist_Y(b,a)$. 

There is a simple relation between the distances in a graph and the distances in its line digraph.  Recall that for an undirected graph, we view each edge as a pair of opposite arcs.

\begin{lemma}\label{lem_dist}
 Let $X$ be a  digraph, and let $(a,b)$ and $(c,d)$ be two arcs of $X$. Then
\[\dist_{\LD(X)}((a,b),(c,d)) =\begin{cases}
\dist_X(c,b)-1,\quad \text{if }(a,b)=(c,d)\\
\dist_X(b,c)+1,\quad \text{otherwise}.
\end{cases}\]
\end{lemma}
\proof
The case where $(a,b)$ coincides with $(c,d)$ is trivial. Suppose $(a,b)$ is different from $(c,d)$. If there is a dipath $u_0=b, u_1, u_2, \cdots, u_m=c$ in $X$ from $b$ to $c$, then $(a, u_0), (u_0, u_1), (u_1,u_2),\cdots, (u_m, d)$ is a dipath in $\LD(X)$ from $(a,b)$ to $(c,d)$. This shows
\[\dist_{\LD(X)}((a, b), (c, d)) \le \dist_X(b,c)+1\]
Now, suppose  $(a, u_0), (u_0, u_1), (u_1,u_2),\cdots, (u_m, d)$ is a shortest dipath in $\LD(X)$ from $(a,b)$ to $(c,d)$. Then for any $i\ne j$, we have $u_i\ne u_j$, as otherwise we could have removed the portion from $(u_i, u_{i+1})$ to $(u_{j-1}, u_j)$ to get a shorter dipath. Hence, $u_0,u_1,u_2,\cdots,u_m$ is a dipath in $X$ from $b$ to $c$, and so
\[\dist_X(b,c)\ge \dist_{LD(X)}((a,b),(c,d))-1.\tag*{\sqr53}\]

\section{Distance regular graphs and their line digraphs}
A graph is \textsl{distance regular} if for any two vertices $u$ and $v$ at distance $\ell$, the number of vertices at distance $i$ from $u$ and distance $j$ from $v$ is a constant $p_{ij}^{\ell}$, which depends only on the distances $i$, $j$ and $\ell$. The parameters $p_{ij}^{\ell}$ are called the \textsl{intersection numbers}.

Suppose our graph $X$ is distance regular of diameter $d$. Let $A_i=A(X_i)$ denote the adjacency matrix of its $i$-th distance graph. It is well-known that
\[\{\seq{A}{0}{1}{d}\}\]
forms an association scheme (see, for example, \cite[Ch 12]{Godsil1993}). In particular,  there are scalars $p_r(i)$, called the \textsl{eigenvalues} of the scheme, such that
\[A_r = p_r(0) A_0 + p_r(1) A_1 + \cdots + p_r(d) A_d,\]
and scalars $q_r(i)$, called the \textsl{dual eigenvalues} of the scheme, such that
\begin{equation}\label{eq:q_r()}
E_r = \frac{1}{|V(X)|} \left(q_r(0) A_0 + q_r(1) A_1 +\cdots + q_r(d)A_d \right).
\end{equation}
Moreover, for any $i$ and $j$, the product $A_iA_j$ lies in the span of the scheme:
\[A_i A_j = \sum_{\ell=0}^d p_{ij}^{\ell} A_{\ell},\]
and $p_{ij}^k$'s are called the \textsl{intersection numbers}.

We can generalize distance regularity to digraphs. A digraph $Y$ is \textsl{distance regular} if for any two vertices $u$ and $v$ with
$\dist_Y(u,v)=\ell$, the number of vertices $w$ such that $\dist_Y(u,w)=i$ and $\dist_Y(w,v)=j$ depends only on the distances $i$,$j$ and $\ell$, and we denote this number by $m_{ij}^{\ell}$.

Now, we are ready to show that the line digraph of a distance regular graph is distance regular.

\begin{lemma}\label{lem_drdg}
Let $X$ be a distance regular graph. Let $Y=\LD(X)$ be its line digraph. Then $Y$ is a distance regular digraph. Moreover, for any $i$ and $j$, we have $m_{ij}^{\ell} = m_{ji}^{\ell}$. 
\end{lemma}
\proof
Let $(a,b)$ and $(c,d)$ be two arcs of $X$. We count arcs $(x,y)$ that are at distance $i$ from $(a,b)$ and from which $(c,d)$ is at distance $j$.  

We will prove the case where $i,j\ge 1$; the other cases follow from a similar argument. By Lemma \ref{lem_dist}, we are counting elements in the set
\[\Phi:=\{(x,y): \dist_X(b,x)=i-1, \dist_X(y, c)=j-1, \dist_X(x,y)=1\}.\]
On the other hand,
\[(A_{i-1} A_1 A_{j-1})_{bc} = \sum_x \sum_y (A_{i-1})_{bx} (A_1)_{xy} (A_{j-1})_{yc}=\abs{\Phi}.\]
Since $X$ is distance regular, there are intersection numbers $p_{1,j-1}^{m}$ and  $p_{i-1,m}^{\ell}$ such that
\begin{align*}
A_{i-1} A_1 A_{j-1} &= A_{i-1} \sum_{m} p_{1,j-1}^{m} A_{m}\\
&=\sum_{\ell} \sum_{m}  p_{i-1,m}^{\ell} p_{1,j-1}^{m} A_{\ell}.
\end{align*}
Therefore, the size of $\Phi$ is a constant $m_{ij}^{\ell}$ that depends only on $i$, $j$ and $\ell$. By Lemma \ref{lem_dist}, $\ell$ is determined by the distance from $(a,b)$ to $(c,d)$. Finally, since $A_{i-1}$ commutes with $A_{j-1}$, we have $m_{ij}^{\ell} = m_{ji}^{\ell}$. 
\qed

Let $Y$ be a digraph. The $i$-th \textsl{distance digraph} of $Y$, denoted $Y_i$, has the same vertex set as $Y$, and vertex $a$ is adjacent to vertex $b$ in $Y_i$ if $\dist_Y(a,b)=i$.  One consequence of Lemma \ref{lem_drdg} is that the  $(0,1)$-adjacency matrices of the distance digraphs of $\LD(X)$ commute, if $X$ is distance regular. Moreover, we have a formula for these adjacency matrices.

\begin{lemma}\label{lem_form}
Let $X$ be a distance regular graph. Let $A_i$ be the adjacency matrix of its $i$-th distance graph. Let $Y=\LD(X)$. Let $Y_i$ be the $i$-th distance digraph of $Y$. Then $A(Y_0)=I$, $A(Y_2) = D_h^T A_1 D_t-I$, and for $i= 1,3,\cdots,d+1$, $A(Y_i) = D_h^T A_{i-1} D_t $. 
\end{lemma}
\proof
Clearly, $A(Y_0)=I$. The remaining statement is a direct translation of Lemma \ref{lem_dist}. Indeed, for $i=1,2\cdots,d$,
\begin{align*}
(D_h^T A_i D_t)_{(a,b), (c,d)} &= e_{(a,b)}^T D_h^T A_i D_t e_{(c,d)}\\
&=e_b^TA_i e_c\\
&=\begin{cases}
1,\quad \text{if }\dist_X(a,b)=i.\\
0,\quad \mathrm{otherwise}.
\end{cases},
\end{align*}
First assume $i=2,3,\cdots, d+1$. We see that $\dist_X(b,c)=i$ if and only if $(a,b)\ne (c,d)$ and 
\[\dist_{\LD(X)}((a,b),(c,d)) = \dist_X(b,c)+1 = i+1,\]
from which it follows that $A(Y_{i+1}) = D_h^T A_i D_t$. Now assume $i=1$. Then $\dist_X(b,c)=i$ if and only if $b$ and $c$ are adjacent. If $(a,b)\ne (c,d)$, then 
\[\dist_{\LD(X)}((a,b),(c,d)) = \dist_X(b,c)+1 = i+1=2;\]
otherwise, 
\[\dist_{\LD(X)}((a,b),(c,d)) = \dist_X(b,c)-1 = i-1=0.\]
Hence $A(Y_2) + I = D_h^T A_1 D_t$.
\qed

Using Lemma \ref{lem_form}, we can verify two properties of $A(Y_i)$, one of which we have already seen.
\begin{corollary}\label{cor_props}
Let $X$ be a distance regular graph. Let $Y=\LD(X)$. Let $Y_i$ be the $i$-th distance digraph of $Y$. Then the following statements hold.
\begin{enumerate}[(i)]
\item $A(Y_0), A(Y_1), \cdots,A(Y_{d+1})$ sum to $J$.
\item $Y$ is distance regular.
\end{enumerate}
\end{corollary}
\proof
We prove (i); the other statement follows similarly. By Lemma \ref{lem_form}, 
\[\sum_{i=0}^{d+1}A(Y_i) = I + (D_h^TA_0D_t) + (D_h^TA_1D_t-I) + \sum_{i=2}^dD_h^TA_{i-1}D_t\]
which reduces to $D_h^TJ D_t$. Since each column of $D_h$ or $D_t$ has exactly one non-zero entry, we have $D_h^TJD_t=J$.
\qed

\begin{figure}[h]
	\centering
\begin{minipage}{0.5\textwidth}
\includegraphics[width=6cm]{LD}
\end{minipage}%
\begin{minipage}{0.5\textwidth}
\includegraphics[width=8cm]{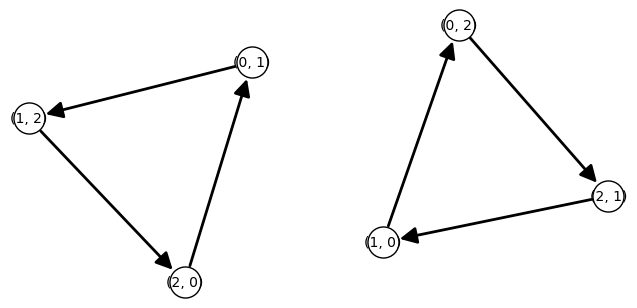}
\end{minipage}
\caption{$\LD(K_3)$ and the underlying oriented graph of $S(\LD(K_3))$}
\label{fig_supp}
\end{figure}

Given a digraph $Y$, let $S(Y)$ denote its \textsl{skew-adjacency matrix}, that is,
\[S(Y) = A(Y) - A(Y)^T.\]
This represents an \textsl{oriented graph}, that is, a digraph with no symmetric pair of arcs.
Note that $A(Y)$ and $A(Y)^T$ may not have disjoint supports, so the underlying oriented graph of $S(Y)$ may differ from $Y$. We illustrate this situation in Figure \ref{fig_supp}.

Two properties of $S(Y_i)$ follow from Lemma \ref{lem_form} and Corollary \ref{cor_props}.
\begin{corollary}\label{cor_SY}
Let $X$ be a $k$-regular distance regular graph. Let $Y_i$ be the $i$-th distance digraph of $\LD(X)$. The following statements hold.
\begin{enumerate}[(i)]
\item $S(Y_1),S(Y_2),\cdots,S(Y_{d+1})$ sum to zero; in particular, they are linearly dependent.
\item For  $i,j=1,2,\cdots,d+1$,
\[S(Y_i) S(Y_j) = S(Y_j) S(Y_i).\]
\end{enumerate}
\end{corollary}
 \proof
We have $S(Y_0)=A(Y_0)-A(Y_0)^T = 0$ and so
\begin{align*}
S(Y_1) + S(Y_2) + \cdots + S(Y_{d+1}) & = \sum_{i=0}^{d+1} (A(Y_i) -A(Y_i)^T)\\
&= J-J^T\\
&= 0.
\end{align*}
Also, for $i=1,2,\cdots, d+1$, 
\[S(Y_i) = D_h^T A_{i-1}D_t - D_t^T A_{i-1} D_h.\]
Since 
\begin{align*}
A_{i-1}D_t D_h^T A_{j-1} &= A_{i-1}D_h D_t^TA_{j-1} = A_{i-1} A_1 A_{j-1} = A_{j-1} A_1 A_{i-1},\\
A_{i-1}D_t D_t^T A_{j-1} &= A_{i-1}D_h D_h^TA_{j-1} =k A_{i-1} A_1 A_{j-1} =k A_{j-1} A_1 A_{i-1},
\end{align*}
it follows that $S(Y_i)$ commutes with $S(Y_j)$.
\qed

\section{Spectral decomposition}
Let $X$ be a $k$-regular distance regular graph of diameter $d$.  Recall that the transition matrix $U$ of the Grover walk on $X$ is
\[
U = \frac{2}{k}A(\LD(X)) - R=R\left(\frac{2}{k}D_t^TD_t-I\right),
\]
In this section, we derive some properties of the eigenprojections of $U$. 

We first cite a result from \cite{Zhan2019}.

\begin{theorem}\label{thm_spd}
Every non-real eigenvalue of $U$ comes from an eigenvalue of $A(X)$. More specifically, let $\lambda$ be an eigenvalue of $A(X)$ that is neither $k$ nor $-k$. Let $E_{\lambda}$ be the orthogonal projection onto the $\lambda$-eigenspace of $X$. Suppose $\lambda=k\cos(\theta)$ for some real number $\theta$. Then $e^{i\theta}$ is an eigenvalue of $U$ with eigenprojection
\[F_{\theta_+} = \frac{1}{2k \sin^2(\theta)} (D_t - e^{i\theta} D_h)^T E_{\lambda} (D_t - e^{-i\theta} D_h),\]
and $e^{-i\theta}$ is an eigenvalue of $U$ with eigenprojection
\[F_{\theta_-}=\frac{1}{2k\sin^2(\theta)} (D_t - e^{-i\theta} D_h)^T E_{\lambda} (D_t - e^{i\theta} D_h).\tag*{\sqr53}\]
\end{theorem}

As a consequence, for any non-real eigenvalue $e^{i\theta}$ of $U$, the difference $F_{\theta_+} - F_{\theta_-}$ is an imaginary scalar multiple of some real skew-symmetric matrix. 

\begin{corollary}\label{cor_Fdiff}
Let $\lambda$ be an eigenvalue of $A(X)$ that is neither $k$ nor $-k$. Let $E_{\lambda}$ be the orthogonal projection onto the $\lambda$-eigenspace of $X$. Suppose $\lambda = k\cos(\theta)$ for some real number $\theta$. Let $F_{\theta_+}$ and $F_{\theta_-}$ be the $e^{i\theta}$-eigenprojection and the $e^{-i\theta}$-eigenprojection, respectively. Then 
\[F_{\theta_+} - F_{\theta_-}=\frac{i}{k\sin(\theta)}(D_t^TE_{\lambda} D_h - D_h^T E_{\lambda} D_t). \tag*{\sqr53}\]
\end{corollary}

We now fix a $\lambda$ that satisfies the condition in Corollary \ref{cor_Fdiff}, and set
\[S_{\lambda} := D_t^TE_{\lambda} D_h - D_h^T E_{\lambda} D_t.\]
We show that if $X$ is distance regular, then $S_{\lambda}$ is a linear combination of the skew adjacency matrices of the distance digraphs of $\LD(X)$ , where the coefficients are differences of dual eigenvalues $q_r(i)$ defined in Equation \eqref{eq:q_r()}.

\begin{lemma}\label{lem_span}
Let $X$ be a distance regular graph on $n$ vertices of diameter $d$. Let $Y = \LD(X)$. Let $Y_i$ be the $i$-th distance digraph of $Y$. If $\lambda$ is the $r$-th eigenvalue of the corresponding association scheme, and $q_r(i)$'s are the dual eigenvalues, then
\[S_{\lambda} = \frac{1}{n}\sum_{i=1}^d (q_r(d) - q_r(i-1)) S(Y_i).\]
%\[S_{\lambda} \in \mathrm{span} \{S(Y_1),S(Y_2),\cdots,S(Y_d)\}\]
\end{lemma}
\proof
Recall that 
\[E_r = \frac{1}{n} \left(q_r(0) A_0 + q_r(1) A_1 +\cdots + q_r(d)A_d \right).\]
If $\lambda$ is the $r$-th eigenvalue in the association scheme, then $E_{\lambda}=E_r$ and
\begin{align*}
S_{\lambda} &= D_t^TE_r D_h - D_h^T E_rD_t\\
&=\frac{1}{n}\sum_{j=0}^d q_r(j) (D_t^T A_j D_h - D_h^T A_j D_t)\\
&=-\frac{1}{n}\sum_{j=0}^d q_r(j) S(Y_{j+1})\\
&=-\frac{1}{n}\sum_{i=1}^{d+1} q_r(i-1) S(Y_i)\\
&=\frac{1}{n}\sum_{i=1}^d (q_r(d)-q_r(i-1)) S(Y_i),
\end{align*}
where the last equality follows from the fact that $S(Y_1), S(Y_2), \cdots, S(Y_{d+1})$ sum to zero.
\qed

\section{Factoring $U^2$}

Every unitary matrix $U$ can be written as the exponential of some skew-Hermitian matrix: if 
\[U = \sum_r e^{i\theta_r} F_r\]
is the spectral decomposition of $U$, then 
\[U = \exp(iH),\]
where $H$ is a Hermitian matrix given by
\[H = \sum_r \theta_r F_r.\]

Now let $X$ be a $k$-regular graph, and $U$ the transition matrix of the Grover walk on $X$. The eigenvalues of $U$ are $\pm1$ and some complex numbers $e^{\pm i\theta_r}$ that come in conjugate pairs. We may assume, without loss of generality, that each $\theta_r$ lies in $[0, \pi]$. Then there is a natural way to group the terms in $H$:
\[H = 0\cdot F_0 + \pi\cdot  F_{\pi} + \sum_{r: 0<\theta_r <\pi}\theta_r \cdot (F_r - \comp{F_r}),\]
where $F_0$ is the $1$-eigenprojection, $F_{\pi}$ is the $(-1)$-eigenprojection, and $F_r$ is the $e^{i\theta_r}$-eigenprojection of $U$, respectively.

From Theorem \ref{thm_spd} and Theorem \ref{cor_Fdiff}  we see that, if $\theta_r\in(0,\pi)$, then $k\cos (\theta)$ is an eigenvalue of $A(X)$, and $F_r - \comp{F_r}$ can be constructed from some eigenspace of $A(X)$. If $-1$ were not an eigenvalue of $U$, then $H$ would take a simpler form
\[H = \sum_{r: 0<\theta_r<\pi} \theta_r (F_r - \comp{F_r}),\]
which likely bears some combinatorial meaning. Unfortunately, the following result from \cite{Zhan2019} shows that $-1$ is always an eigenvalue for $U$.

\begin{lemma}
Let $X$ be a $k$-regular graph on $n$ vertices. Then $(-1)$ is an eigenvalue of $U$, with multiplicity $nk/2-n+2$ if $X$ is bipartite, and $nk/2-n$ otherwise.
\qed
\end{lemma}

The square of $U$ however may not have $-1$ as an eigenvalue. In fact, if $A(X)$ is invertible, then $\pm i$ are not eigenvalues of $A(X)$, and so $-1$ is not an eigenvalue of $U^2$. In this case, we have
\[U^2=  \exp(iH),\]
where 
\[H = 2 \sum_{r: 0<\theta_r<\pi} \theta_r (F_r- \comp{F_r}).\]

Our discussion in the last section shows that $H$ is a linear combination of the skew-adjacency matrices of the distance digraphs of $\LD(X)$.

\begin{theorem}
Let $X$ be a distance regular graph with diameter $d$. Let $U$ be the transition matrix of the Grover walk on $X$. Let $Y_j$ be the $j$-th distance digraph of $\LD(X)$. If $A(X)$ is invertible, then there are real scalars $\seq{t}{1}{2}{d}$ such that
\[U^2 = \exp(t_1S(Y_1)) \exp(t_1 S(Y_2))\cdots \exp(t_d S(Y_d)).\]
 Moreover, if $n$ is the number of vertices in $X$, and $p_r(j)$'s and $q_r(j)$'s are the eigenvalues and dual eigenvalues of the association scheme, then
\[t_j = \frac{2}{n} \sum_{r: -k<p_1(r)<k} \frac{\arccos(p_1(r)/k)}{\sqrt{k^2-(p_1(r))^2}}( q_r(j-1)-q_r(d)).\]
\end{theorem}
\proof
We have
\[U^2 = \exp(iH)\]
where 
\[H = 2 \sum_{r: 0<\theta_r<\pi} \theta_r (F_r- \comp{F_r}).\]
 By Lemma \ref{lem_span}, for each eigenvalue $p_1(r)\in(-k, k)$, 
\[S_{p_1(r)} = \frac{1}{n} \sum_{j=1}^d (q_r(d) - q_r(j-1)) S(Y_j),\]
and so by Corollary \ref{cor_Fdiff},
\begin{align*}
H &= 2i \sum_{r: -k < p_1(r) < k} \frac{\arccos(p_1(r)/k)}{\sqrt{k^2-(p_1(r))^2}} (D_t^T E_{p_1(r)} D_h - D_h^T E_{p_1(r)} D_t)\\
&=\frac{2i}{n} \sum_{r: -k<p_1(r)<k} \frac{\arccos(p_1(r)/k)}{\sqrt{k^2-(p_1(r))^2}}\sum_{j=1}^d (q_r(d)-q_r(j-1)) S(Y_j).
\end{align*}
Finally, by Corollary \ref{cor_SY}, these skew-adjacency matrices $S(Y_j)$ commute, so we can factor $U^2$ into the desired form.
\qed

If we rewrite $\exp(t_iS(Y))$ in the above theorem as $\exp(i t_i (iS(Y_i)))$, then we have the transition matrix of a  continuous-time quantum walk on an oriented graph, evaluated at time $t_i$. Thus, for a distance regular graph with diameter $d$ and invertible adjacency matrix, the Grover walk on it factors into at most $d$  continuous-time quantum walks on the distance digraphs of $\LD(X)$. In particular, a  discrete-time quantum walk on $K_n$ is equivalent to a  continuous-time quantum walk on its line digraph.

\begin{corollary}
Let $U$ be the Grover walk on $K_n$. Then 
\[U^2 = \exp(tS(\LD(K_n)))\]
where
 \[t = \frac{2\arccos(1/(1-n))}{\sqrt{n^2-2n}}.\tag*{\sqr53}\]
\end{corollary}
 \proof
The association scheme for $K_n$ is 
\[\{A_0=I, \quad A_1= J-I\},\]
with eigenprojections
\[E_0=\frac{1}{n}J, \quad E_1=I -\frac{1}{n}J,\]
eigenvalues
\[p_0(0)=p_0(1)=1, \quad p_1(0)=n-1,\quad p_1(1)=-1,\]
and dual eigenvalues
\[q_1(0)=q_0(1)=1,\quad q_1(0)=n-1,\quad q_1(1)=-1.\]
The result then follows from a direct calculation.
\qed

\begin{figure}[h]
	\centering
\begin{minipage}{0.5\textwidth}
\includegraphics[width=6cm]{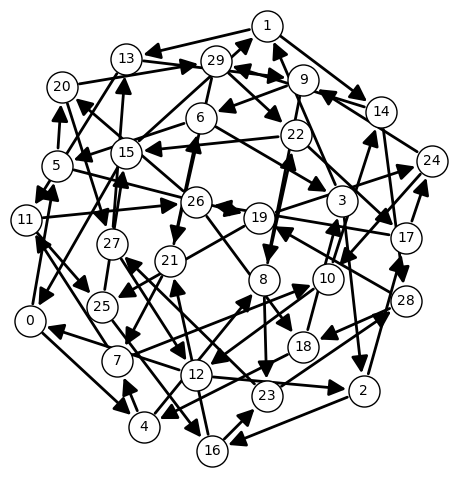}
\end{minipage}%
\begin{minipage}{0.5\textwidth}
\includegraphics[width=6cm]{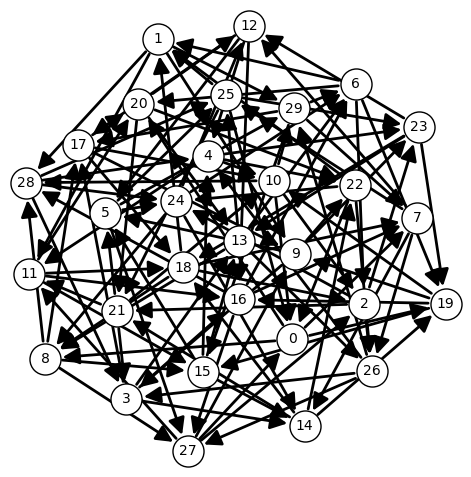}
\end{minipage}
\caption{Underlying oriented graphs of $S(\LD(P))$ and $S(\LD(P)_2)$}
\label{fig_LDPet}
\end{figure}

The Grover walk on an invertible strongly regular graph factors into at most two  continuous-time quantum walks. For the Petersen graph $P$, the factors of $U^2$ are  continuous-time quantum walks on the oriented graphs shown in Figure \ref{fig_LDPet}.

Finally, we remark here that not all $d$ distance digraphs may show up in the factorization. For example, if $U$ is the Grover walk of the $3$-cube $Q_3$, and $Y_2$ is the second distance digraph of $\LD(Q_3)$, then there is a scalar $t$ such that $U^2 = \exp(t S(Y_2))$.   It will be interesting to characterize distance regular graphs on which the discrete-time quantum walk factors into fewer than $d$ continuous-time quantum walks.

\section{Implementing the search algorithm}
We show how the aforementioned factorization allows us to implement some search algorithms using continuous-time quantum walks. To start, we introduce a quantum walk based search algorithm due to Shenvi, Kempe and Whaley \cite{Shenvi2003}. Given a graph $X$ and a marked vertex $a$, the \textsl{oracle}, denoted $O_a$, is the diagonal matrix indexed by the arcs of $X$, such that
\[(O_a)_{(u,v), (u,v)}= 
\begin{cases}
-1,&\quad u=a,\\
1,&\quad u\ne 1.
\end{cases}\]
Assume $X$ is $k$-regular on $n$ vertices. The state space associated with $X$ is $\cx^n\otimes \cx^k$, the vector space of complex-valued functions on the arcs, and the oracle $O_a$ can be written as
\[O_a = (I-2E_{aa}) \otimes I.\]
To find the marked vertex $a$, we apply the oracle $O_a$ and the transition matrix $U$ of the Grover walk alternately to the following initial state:
\[x_0 = \frac{1}{\sqrt{nk}} \one_n \otimes \one_k.\]
Thus, at time $t$, the system is in state
\[x_t = (U O_a)^t x_0.\]
If we measure at time $t$, the outcome will be an arc $(u,v)$ with probability
\[\abs{\ip{x_t}{e_{(u,v)}}}^2.\]
We call the sum
\[\sum_{v\sim a} \abs{\ip{x_t}{e_{(a,v)}}}^2\]
the \textsl{success probability} of finding the marked vertex at time $t$.

The well-known Grover's search algorithm \cite{Grover1996} is equivalent to finding a marked vertex using the above procedure - the underlying graph is a complete graph $K_n$ with a loop at each vertex, and the success probability is reasonably high at time $t=O(\sqrt{n})$. 

While one walk step per oracle query seems standard in search algorithms, Wong and Ambainis \cite{Wong2015} proposed a modification, where the algorithm takes multiple walk steps per oracle query. That is, the transition matrix in this algorithm is 
\[U^k O_a\]
for some positive integer $k$. Thus, if $k$ is even and $X$ is an invertible distance regular graph, we can implement the walk steps $U^k$ using continuous-time quantum walks. As for the oracle query, note that $O_a$ has spectral decomposition
\[O_a = 1\cdot ((I-E_{aa}) \otimes I) + (-1)\cdot (E_{aa} \otimes I),\]
and so 
\[O_a = \exp(i\pi (E_{aa} \otimes I)).\]
The underlying graph of $E_{aa}\otimes I$ consists of isolated vertices, each representing an arc of $X$, and every vertex representing an outgoing arc of $a$ has a loop on it. Thus, we can implement $O_a$ by running a continuous-time quantum walk on this graph for the time period of $\pi$.

\section{Generalization to association schemes, and future work}
We have seen that for an invertible distance regular graph $X$, a  discrete-time quantum walk on $X$  factors into  continuous-time quantum walks on the distance digraphs of the line digraph of $X$. It will be interesting to translate properties on the discrete side, such as perfect state transfer, into properties on the continuous side, and vice versa.

In general, if the adjacency matrix of a graph $X$ is invertible and lies in the Bose Mesner algebra of some association scheme, then we have a similar factorization.

\begin{theorem}
Let $\cA=\{\seq{A}{0}{1}{d}\}$ be an association scheme. Let $D_t$ and $D_h$ be the arc-tail and arc-head incidence matrices of $X_1$. For each $i=0,1,\cdots, d$, define a skew-symmetric  $(1, -1)$-matrix by
\[S_i = D_h^T A_i D_t - D_t^T A_i D_h.\]
Then $S_iS_j = S_j S_i$.

Let $X$ be a graph whose adjacency is invertible and lies in the Bose Mesner algebra of $\cA$. Let $U$ be the transition matrix of the Grover walk on $X$. Then there are real scalars $\seq{t}{1}{2}{d}$ such that
\[U^2 = \exp(t_1 S_1) \exp(t_2 S_2) \cdots \exp(t_d S_d).\]
\end{theorem}
\proof
The fact that $S_i$ commutes with $S_j$ follows from Lemma \ref{lem_props} and that $\seq{A}{0}{1}{d}$ commute with each other.

For the second part, note that the Bose Mesner algebra of $\cA$ has a basis $\{\seq{E}{0}{1}{d}\}$ of idempotents. By the spectral decomposition of $U$ and Corollary \ref{cor_Fdiff}, we see that $U^2 = \exp(iH)$, where $H$ lies in
\[\mathrm{span}\{D_t^T E_i D_h - D_h^T E_i D_t: i=0,1,\cdots, d\},\]
that is,
\[\mathrm{span}\{D_t^T A_i D_h - D_h^T A_i D_t: i=0,1,\cdots, d\}.\tag*{\sqr53}\]

When $X$ is distance regular, the matrix $S_i$ in the above theorem has a nice interpretation: it is the skew-adjacency matrix of the $i$-th distance digraph of the line digraph of $X$. In the more general situation, however, the relation between $S_i$ and $X$ has not been explored. It will be helpful to look into a few common association schemes, such as groups schemes, and study the ``lift" of the Schur idempotents via the the arc-tail and arc-head incidence matrices.

\section{Acknowledgment}
The author would like to thank Ada Chan for her inspiring discussion, and David Roberson for bringing up the generalization to association schemes. The author acknowledges the support of the York Science Fellow program. The author thanks the anonymous referees for their valuable comments.

\bibliographystyle{amsplain}

\providecommand{\bysame}{\leavevmode\hbox to3em{\hrulefill}\thinspace}
\providecommand{\MR}{\relax\ifhmode\unskip\space\fi MR }

\providecommand{\MRhref}[2]{%
  \href{http://www.ams.org/mathscinet-getitem?mr=#1}{#2}
}
\providecommand{\href}[2]{#2}

\end{document}